\documentclass[12pt]{amsart}

\usepackage[latin1]{inputenc}

\usepackage{amsmath}
\usepackage{amssymb}
\usepackage{graphicx}
\usepackage{graphics}
\usepackage{amsthm}
\usepackage{amscd}
\usepackage{color}
\usepackage{amsfonts}
\usepackage{fancyhdr}
\newtheorem{theorem}{Theorem}[section]

\newtheorem{mainthm}{Theorem}
\newtheorem*{theorem*}{Theorem}
\newtheorem{corollary}[theorem]{Corollary}

\newtheorem{lemma}[theorem]{Lemma}
\newtheorem{Remark}[theorem]{Remark}
\newtheorem*{definition*}{Definition}
\newtheorem*{af}{Assertion}
\newtheorem{definition}[theorem]{Definition}
\newtheorem*{Question*}{Question}
\newtheorem{Question}{Question}

\def\N{\mathbb{N}}
\def\Z{\mathbb{Z}}

\def\R{\mathbb{R}}

\def\tq{\,\mid\,}
\def\set#1{\left\{\, #1 \,\right\}}

\def\norm #1{\Vert \,#1\, \Vert\,}
\makeatletter

\newcommand{\Rmnum}[1]{\expandafter\@slowromancap\romannumeral #1@}
\makeatother

\begin{document}

\vspace{-2cm}
\title[New examples]{Hyperbolic sets that are not   contained \\in a locally maximal one}

\author{Adriana da Luz}

\vspace{-2cm}
\maketitle

\begin{abstract}
In this paper we study two properties related to the structure of hyperbolic sets.
First we construct new  examples answering in the negative the following question
posed by  Katok and Hasselblatt in [\cite{hk}, p. 272]
\begin{Question*} Let $\Lambda$ be a hyperbolic set, and let $V$ be an open neighborhood of
 $\Lambda$. Does there exist a locally maximal hyperbolic set  $\widetilde{\Lambda}$
such that $\Lambda \subset \widetilde{\Lambda} \subset V  $?
\end{Question*}
We show that such examples are present in linear anosov diffeomorophisms of  $\mathbb{T}^3$,
and are therefore robust.

Also we construct  new examples of  sets that are not contained in any locally maximal hyperbolic set.
The examples known until now were constructed by Crovisier in \cite{c} and by Fisher in \cite{Fi}, and these were either in dimension bigger than 4 or they were not transitive.
We give a transitive and robust example in
 $\mathbb{T}^3$. And show that such examples cannot be build in dimension 2.

\end{abstract}

\let\thefootnote\relax\footnote{The autor was partially supported by CSIC.}

\section{\textbf{Introduction}}\label{sec:1}

In the '60s, Anosov (\cite{A2}) and Smale (\cite{S}) began the study of some compact  invariant sets,
 whose tangent space splits into invariant,
uniformly contracting and uniformly expanding directions. 
More precisely, a hyperbolic set is defined to be a compact invariant subset of a  compact manifold $\Lambda\subset M$
of a diffeomorphism $f$ such that the tangent space at every $x\in\Lambda$ admits an invariant
 splitting that satisfies:
\begin{itemize}
  \item $T_x M=E^s(x)\oplus E^{u}(x)$
   \item $Df_x(E^s(x))=E^s(f(x))$ y  $Df_x(E^{u}(x)=E^{u}(f(x))$
  \item  there are constants $C>0$ and $\lambda\in (0,1)$
   such that for every $n\in\N$ one has $\norm{Df^n(v)}\leq C\lambda^n\norm{v}$ for $v\in E^s(x)$  and
    $\norm{Df^{-n}(v)}\leq C\lambda^{-n}\norm{v}$ for $v\in E^u(x)$ .

\end{itemize}

A specially interesting case is when the
hyperbolic set $\Lambda$ is the non-wandering set of $f$.
Particularly when we also have that the set of periodic points of $f$
is dense in the non-wandering set $\Omega(f)$, we say that $f$ is Axiom A.



Given the relevance of these diffeomorphisms in the study of hyperbolic dynamics is natural to ask
 what kind of sets may or may not be a basic pieces of some spectral decomposition.

All basic pieces have the following property:

\begin{definition}\label{locprod}
Let $ f: M \to M $ be a diffeomorphism and $ \Delta $ a compact invariant hyperbolic set.
 We say that $ \Delta $ has locald product structure exists $ \delta> 0 $ such that if
 $ x, y \in \Delta $, and $ d (x, y) <\delta $ then
 $ W^s_{\varepsilon} (x) \cap W^u_{\varepsilon} (y) \in \Delta $ where
$ \varepsilon $ is as in the stable manifold theorem.
\end{definition}
We will focus now on whether or not a set has this property. If they do not, we will be interested
 in studying  whether or not the set  is contained in an other set  having this property.

Sets having this property are interesting on themselves since they can be thought of,
locally, in coordinates of the stable and unstable manifold of a point.
Also this property is equivalent to others that are very useful to understand the dynamics
 of a neighborhood of the set. Some of them are having the shadowing property or being locally maximal.

 Many of the best known examples of hyperbolic compact sets have this property.
 Some  examples could be the solenoid, the torus under an Anosov diffeomorphism or a horseshoe.
Also, there are examples of simple sets that
 do not verify this property, for instance,  the closure of the orbit of a homoclinic point.
 However, for a long time all the  examples known of sets that did not have
 local product structure could be included in a set having this property.
 Moreover all known examples had such a set included in any neighborhood of the original one.
 In the 1960's  Alexseyev  asked the following question
 (that was later posed by  Katok and Hasselblatt in [\cite{hk}, p. 272])
\begin{Question}\label{q1} Let $\Lambda$ be a hyperbolic set, and let $V$ be any open neighborhood of
 $\Lambda$. Does there exist a locally maximal hyperbolic set  $\widetilde{\Lambda}$
such that $\Lambda \subset \widetilde{\Lambda} \subset V  $?
\end{Question}

Also the following related question was  unanswered:
\begin{Question}\label{q2}
Given a hyperbolic set $\Lambda$ does there exist a hyperbolic compact invariant
set with local product structure
such that $\Lambda\subset\widetilde{\Lambda} $?
\end{Question}

Both questions remained open  until 2001, when Crovisier \cite{c} constructed an example based
on an example of Shub in \cite{hps} that answer question \ref{q2} in the negative
(and therefore question \ref{q1}). This example is on the 4-torus.

Later, Fisher \cite{Fi} constructed several other examples of this sort.
He constructed robust examples in any dimension, and transitive volume preserving examples in dimension 4.

In spite of this there are still some natural questions left to answer
\begin{itemize}
  \item  Does there exist an open set $U$ (in the $C^1$ topology) of diffeomorphisms
such that every $f \in U$ possesses an invariant transitive hyperbolic set that is not
contained in a locally maximal one   on any manifolds?
\item  Does there exist robust and transitive examples answering Question \ref{q1} in the negative on
 manifolds with dimension lower than 4?

\item Does there exist an example answering Question \ref{q1} in the negative but that it is contained in a bigger  set having local product structure?
\end{itemize}

 In section (3) we will show
 \begin{mainthm}\label{main1}
 Let $f_A:\mathbb{T}^3\to\mathbb{T}^3$ be an Anosov diffeomorphism. There is a connected, compact proper inavriant subset of $\,\mathbb{T}^3$, such that the only locally maximal set containing it is $\,\mathbb{T}^3$.
 \end{mainthm}
This answers our last question.
 Note that the same will be true for any $g$ sufficiently close to $f_A$.
 We also note that constructing this kind of examples is not possible for $\mathbb{T}^2$ since all invariant compact proper sets are 0-dimensional and from \cite{A1}, in any neighborhood there is a locally maximal set that contains them.

 In section (4) we describe  a well known example by Mañe in \cite{M} that we will use on section (5) to construct a new
 example of a set that is not included in any locally maximal set. This example gives a partial answer to our second question.  It is robust, transitive, 
 and it is a 3 dimensional example, which shows there are examples of this in lower dimensions. The previous  examples
had either tangencyes or came from skew-products 
so, they where not transitive or where in $dim\geq4$.
In (5) we proved the following:
\begin{mainthm}\label{main2}
 There exists $\mathcal{U}\subset Diff(\mathbb{T}^3)$ such that
for every $g\in U$ there is a there is a  compact, proper, invariant, hyperbolic subset of $\mathbb{T}^3$, such that there is no  locally maximal set containing it.
\end{mainthm}

In the case of 2 dimensional surfaces our first 2 questions can be combined in the following
\begin{Question*}
If $dim(M)=2$, and $\Lambda\subset M $ is a transitive hyperbolic set and $U$ is any given neighborhood of $\Lambda$. Does there exist compact invariant
set with local product structure
such that $\Lambda\subset\widetilde{\Lambda}\subset U $?
\end{Question*}

We will give a positive answer to this question. In section (6) we will show:

\begin{mainthm}\label{main3}
 Let $f: M\to M$ be a  diffeomorphism, $M$ a compact surface  and $\Lambda\subset M$  a  compact hyperbolic invariant set. If we also have that  $\Omega f\mid_{\Lambda}=\Lambda$ then for any neighborhood  $V$ of $\Lambda$, there exist  $\widetilde{\Lambda}$ such that  $\widetilde{\Lambda}$ is compact hyperbolic invariant and with local product structure and, $$\Lambda\subset\widetilde{\Lambda}\subset V\,.$$

\end{mainthm}


 \section{Preliminaries}
 Let $M$ be a compact manifold, $f$ a $C^r$ diffeomorphism, and $\Lambda$ a hyperbolic set.

 For $\varepsilon > 0$ sufficiently  small and $x \in \Lambda$ the local stable and unstable manifolds
are respectively:
$$W_{\varepsilon}^s(x, f)=\set{y \in M | \text{ for all } n \in \N, d(f^n(x), f^n(y)) < \varepsilon}\,,$$ and
$$W_{\varepsilon}^u(x, f)=\set{y \in M | \text{ for all } n \in \N, d(f^{-n}(x), f^{-n}(y)) < \varepsilon}\,.$$
The stable and unstable manifolds are respectively:
$$W^s(x, f) =\bigcup_{n\geq0}f^{-n}(W_{\varepsilon}^s(f^n(x), f))\,,$$ and $$W^u(x, f) =\bigcup_{n\geq0}f^{n}(W_{\varepsilon}^u(f^{-n}(x), f))\,.$$

The stable and unstable manifolds are $C^r$
injectively immersed submanifolds.
If two points of $\Lambda$ are sufficiently close, The local stable and unstable manifolds
intersect transversely at a single point.

A very useful property of hyperbolic set is the following:
\begin{definition}
 Let $f: M \to M$ be a diffeomorphism and $\alpha> 0$. We say that $\set{x_n}_{ n \in \Z}$ is an
  $\alpha$-pseudo orbit (for $f$)
  if $d (f (x_n), x_{n +1}) \leq\alpha$ for all $n \in \Z$.
\end{definition}
\begin{theorem} (Shadowing Lemma). Let $f:M \to M$ be a diffeomorphism and
$\Lambda$ a compact hyperbolic set. Then, given  $\beta> 0$, there exists $\alpha> 0$
 such that every $\alpha$-pseudo orbits in $\Lambda$ is $\beta$ shadowed by
 an orbits (not necessarily in $\Lambda$).That is, if ${x_n} \in \Lambda$ is a
  $\alpha$-pseudo orbit, then there exists
 $y \in M$ such that $d (f^n   (y), x_n) \leq \beta$ for all $n \in \Z$.
 \end{theorem}
Let us recaall the following  definition:\begin{definition*}\ref{locprod}
Let $ f: M \to M $ be a diffeomorphism and $ \Lambda $ a compact invariant
 hyperbolic set. We say that $ \Lambda $ has local product structure exists
  $ \delta> 0 $ such that if $ x, y \in \Lambda $, and $ d (x, y)  <\delta $
then $ W^s_{\varepsilon} (x) \cap W^u_{\varepsilon} (y) \in \Lambda $.
\end{definition*}
As a consequence of the shadowing theorem we have:
\begin{corollary}
If in addition  to the other hypothesis we have that  $\Lambda $ has local product structure, then every
$\alpha$-pseudo orbits in $\Lambda$ is $\beta$ shadowed by  an orbits  in $\Lambda$.
\end{corollary}

 With this we can show a very important equivalence with having local product structure that
 is being locally maximal :
 \begin{definition}
 A hyperbolic set $\Lambda$ is called locally maximal (or isolated) if there exists a neighborhood
 $V$ of $\Lambda$ in $M$ such that $\Lambda =\bigcap_{n\in\Z}f^n(V)$.
 \end{definition}

 \begin{corollary}
 A hyperbolic set $\Lambda$ is locally maximal if and onely if $\Lambda $ has local product structure.
 \end{corollary}
 As in \cite{A}  we name the properties we are going to be dealing with.
 \begin{definition}
   We say that a hyperbolic set $ \Lambda \subset M$ is premaximal,
   if there exists
a  hyperbolic set $ \Delta \subset M $ with local product structure (maximal invariant set)
such that $ \Lambda \subset \Delta $.
  \end{definition}
  \begin{definition}
We say that a hyperbolic set $ \Lambda $ is locally premaximal, if for every neighborhood $ U $ of
$\Lambda$, there is  a hyperbolic set $ \Delta $ with  local product structure
 such that $ \Lambda \subset \Delta \subset U $.
\end{definition}


 \section{Proof of Theorem \ref{main1}: A set that is not locally premaximal}\label{sec2}
 In this section we prove that there is a subset of the $\mathbb{T}^3$,
 invariant under a linear Anosov diffeomorphism $f$,
 that is not locally premaximal.

Let $f$ be a Anosov diffeomorphism  in $\mathbb{T}^3$ that is induced form
 $ A \in GL(3, \Z)$ which is a hyperbolic toral automorphism with only one eigenvalue grater than one, and
all eigenvalues real, positive, simple, and irrational. Let $\pi:\R^3 \to\mathbb{T}^3$ be such that $\pi \circ A=f\circ\pi$.
Let us also suppose that $f$ has two fixed points $x_0$ and $x_1$, and $\pi(0,0,0)=x_0$.
As a consequence of the results in  \cite{ha} we have:
\begin{theorem}\label{hankok}
Let $f:\mathbb{T}^3\to\mathbb{T}^3$ be a hyperbolic automorphism, we can find a path $\gamma$ in $\mathbb{T}^3$, such that the set $\overline{\mathcal{O}(\gamma)}\subsetneq\mathbb{T}^3$, is compact, connected and non trivial.
\end{theorem}
This curve can
 also be constructed so that it's image contains a fixed point.
For this $\gamma$ we note $\Lambda=\overline{\mathcal{O}(\gamma)}$.

We will prove now that in this conditions the only set with local product structure
containing $\Lambda$ is the whole $\mathbb{T}^3$,  following mainly the ideas in \cite{M2}.
 Here Mañe proves that every compact, connected, locally maximal subset of
$\mathbb{T}^n$ under a linear hyperbolic automorphism must be of the form
$\Delta=x+G$, where $x$ is a fixed point
and $G$ is an invariant compact subgroup.
In particular in dimension $3$ this implies that $\Delta=\mathbb{T}^3$ or $\Delta=x$ .
We will adapt the proof to the case where $\Delta$ is not connected
but contains non trivial compact, connected, invariant set
that  contains a fixed point.

\begin{definition}
Let $ \Lambda \subset \mathbb{T}^3$ be a compact, connected and invariant, such that $ x_0 \in \Lambda $.
 We say that a curve $ \gamma: [0,1] \to \mathbb{T}^3 $ is $ \delta $-adapted to $ \Lambda $,
if there are $ 0 = t_0 <t_1 <\dots <t_m =1 $  such that $ \gamma (t_j) \in \Lambda $,
 and $ d(\gamma (t), \gamma (t_j)) <\delta $
for all $ t_j \leq t \leq t_ {j +1} \, \, $ and $ 0 \leq j \leq m $
\end{definition}

We define $ \Gamma_{\delta} $ as the subgroup of $ \pi_1 (\mathbb{T}^3, e) = \Z ^3 $
generated  by arcs
$ \gamma: [0,1] \to \mathbb {T}^3 $, $ \delta $-adapted such that $ \gamma (0) = \gamma (1) = x_0 $.
 Note that if $ \delta_1 <\delta_2 $ then $ \Gamma_ {\delta_1} \subset \Gamma_ {\delta_2} $.

 Using the continuity of $A $ we have that, given $ \delta $
there is a $ \delta_1 $, such that $ A(\Gamma_ {\delta '}) \subset \Gamma_ {\delta} $ for all
$ 0 <\delta '<\delta_1 $.

The idea now is to define a $ \Gamma_0$ which we would naively define as the subgroup limit of
 $ \Gamma_ {\delta} $ with $ \delta $
 going to zero.  A first attempt to define it would consider
  $ \bigcap_{\delta> 0} \Gamma_ {\delta} $ but
that set might empty and not represent what we want it to.
 Instead we define $ N_{\delta} $ as the subspace
of $ \R^3 $  generated $ \Gamma_{\delta} $.
 We define $ N_0 = \bigcap_{\delta> 0} N_ {\delta} $ and
 $ \Gamma_0 = N_0 \cap \Z^3 $. Note that $ A(N_0) =  N_0$.

\begin{lemma} \label{lema3}
In the above mentioned conditions , $ (N_0 / \Gamma_0) $ is $ \mathbb{T}^3 $ or $ x_0 $.
\end{lemma}

\proof
First we note that $ (N_0 / \Gamma_0) $  is $f$-invariant since $A(N_0)=N_0$, and $N_0\cap\Gamma_0$ so $ (N_0 / \Gamma_0) $ is an invariant sub-torus.
A result from \cite{h} tells us that if the stable or unstable manifold are 1 dimensional then the only  connected, locally connected, compact, invariant hyperbolic subsets are fixed points and the whole torus

\endproof

Note that since $\Gamma_0=\Z^3\cap N_{0}$ then the previous lemma implies that $\Gamma_0=\Z^3$ or $\Gamma_0=0$.

\begin{lemma}\label{remark}
If $\Gamma_0=\Z^3$ and $\pi^s: N_0\to E^u$  is  the  projection  associated  with  the  splitting $N_0 =
E^s\oplus E^u$,  then there exists $\delta_0$ such that $\pi^s(\Gamma_{\delta})$  is  dense  in  $E^u$ for  all  $0 < \delta <\delta_0$.
\end{lemma}
\proof
To see this, note that $N_{\delta_{1}}\subset N_{\delta_{2}}$ if $\delta_1\leq \delta_2$.
This implies that for some $\delta_0$, $\,N_{\delta}=N_{0}$ for all $0<\delta\leq\delta_0$.
It follows that $$\Gamma_{\delta}\subset\Z^3\cap N_{\delta}=\Z^3\cap N_{0}\,.$$
On the other hand $dim(N_0)=dim(N_{\delta})=ran(\Gamma_{\delta})$, so $ran(\Gamma_{\delta})=3$ and there is an isomorphism $\phi:\Gamma_{\delta}\to\Z^3$.

If $a\in\pi^s(\Gamma_{\delta})$, and since $E^u+a$ is irrational, there is a unique $a'\in\Gamma_{\delta}$ such that $\pi^s(a')=a$. If there where $a'$ and $a''$ such that $\pi^s(a')=\pi^s(a'')=a$, then $a''=E^s+a'$. This is impossible since $a',a''\in\Z^3$ and  $E^s+a'$ is a totally irrational plane.

We define now $\varphi:\pi^s(\Gamma_{\delta})\to\pi^s(\Z^3)$ as $\varphi(a)=\pi^s(\phi(a'))$
which is an isomorphism.


\endproof

Now we consider $\Lambda$ to be the set described by Hancock (\ref{hankok}).Then $\Lambda$ is
compact, connected, invariant, it
contains a fixed point $x_0$ and is not trivial.
Let us suppose there exists a set $\Delta$ with local product structure containing $\Lambda$,
 and let us call
it's lift $\widehat{\Delta}$.

The strategy now is to see that  such a $\Delta $, must contain a dense set in the unstable
 manifold of $ x_0$
(which is of dimension 1). Since $\Delta $ is compact then $ \Delta  = \mathbb {T}^3 $.
\begin{definition}
We say that $ x $ and $ y $ are n-$ \varepsilon $-related in $\widehat{\Delta}$ if there exists
 sequence of point $ x = x_0, x_1, \dots, x_n = y $ such that:

\begin{itemize}
\item $ x_i \in \widehat{\Delta} $ for $ i = 1, \dots, n $
\item $ \pi^s(x_ {i +1}-x_ {i}) \leq \varepsilon $ for $ 1 \leq i \leq n $
\item $ \pi^u(x_ {i +1}-x_ {i}) \leq \varepsilon $ for $ 1  \leq i \leq n $
\end{itemize}
\end{definition}

\begin{lemma} \label{EPG}
If $ x, y \in \widehat{\Delta}$ are n-$\varepsilon $-related, with  $ \varepsilon $ sufficiently small,
 then $ (x + E^s) \cap (y + E^{u} ) \in \widehat{\Delta} $.
\end{lemma}
\proof
We take $ \varepsilon <\delta $ with $ \delta $ from the local product structure. We prove this lemma
by induction. For $ n = 1 $ the  property is verified  by the local product structure.
\begin{figure}[h]
\begin{center}
\scalebox{0.35}{\includegraphics{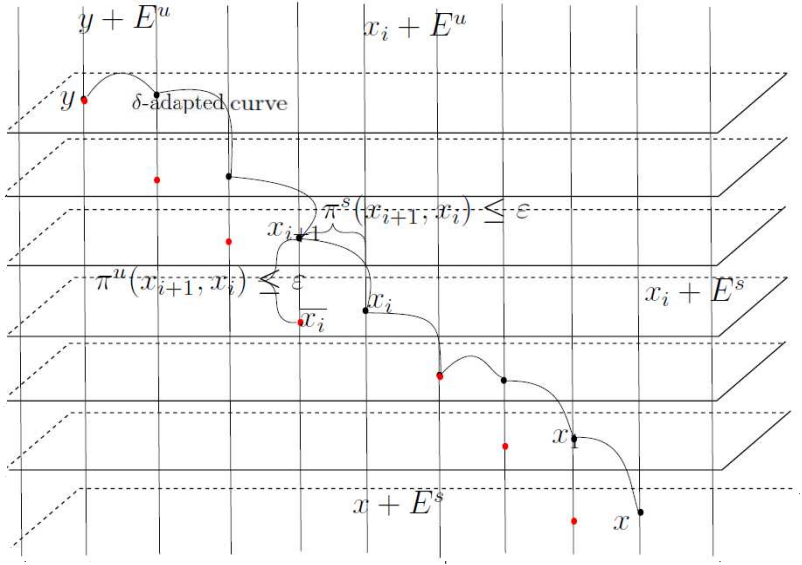}}
\caption{A n-$ \varepsilon $-relation between  $ x $ and $ y $.  }\label{lpgl}
\end{center}
\end{figure}

Suppose now that $ x, y \in \widehat{\Delta} $ are
n-$ \varepsilon $-related. We have $ x = x_0, x_1, \dots, x_n = y $ as in the definition. We define
$ \overline{x_j} = (x_j + E^s) \cap (x_ {j +1} + E^{u}) $ for $ 0 \leq j \leq n-1 $ 
 Note that
$ \overline{ x_0} $ and $ \overline{x_ {n-1}} $ are (n-1)-$\varepsilon $-related because:
\begin{itemize}
\item $ \overline{x_j} \in \widehat{\Delta} $ for $ j = 1, \dots, n-1 $ by induction hypothesis,
\item $ \pi^s (\overline{x_ {j +1}} - \overline{x_ {j}}) = \pi^s (x_ {j +1}-x_ {j}) \leq \varepsilon $ for
 $ 1 \leq j \leq n-1 $,
\item $ \pi^u (\overline{x_ {j +1}} - \overline{x_ {j}}) = \pi^u(x_ {j +1}-x_ {j}) \leq \varepsilon $ for
$ 1 \leq j \leq n-1 $.
\end{itemize}

Then we then have, $ (\overline{x_0} + E ^ s) \cap (\overline{x_ {n-1}} + E^{u}) = z \in \widehat{\Delta} $.

Since we also have $ (\overline{x_0} + E^s) = (x_0 + E^s) $ and
$ (\overline{x_ {n-1}} + E^{u}) = (x_n + E^{u}) $,
we conclude that $ (x_0 + E^s) \cap (x_n + E^{u}) = z \in \widehat{\Delta} $.
\endproof
The following theorem implies theorem \ref{main1}.
\begin{theorem}
Let $ \Lambda $ be a compact, connected, invariant, such that $ x_0 \in \Lambda $, and $ x_0 \neq \Lambda .$
 Suppose there is $ \Delta $ such that $\Lambda\subset\Delta$
 and $\Delta$ is compact invariant and with local product structure. Then $ \Delta =\mathbb{T}^3 $.
\end{theorem}

\proof
Let $ \widehat {\Delta} $ and $ \widehat {\Lambda} $ be the lifts of  $ \Delta $ and $ \Lambda $ respectively.

If $ \Delta $ is compact invariant and local product structure, then by Lemma \ref{EPG},
 if we have two points $ x, y \in \widehat {\Delta} $
 which are  n-$ \varepsilon $-related, we have $ (x + E^s) \cap (y + E^{u}) \in \widehat {\Delta} $.

The goal then is to see that $x_0 $ and any point $ \Gamma_{\delta} $ are n-$\varepsilon $-related and
therefore $ \pi^s ( \Gamma_{\delta}) \subset \widehat {\Delta} $. Since $ \pi^s ( \Gamma_{\delta})$ by \ref{remark} is dense in $E^u$, then $$\overline{ \pi^s (\Gamma_{\delta})}=E^u\subset\widehat {\Delta}\text{ and }\mathbb{T}^3 =\overline{\pi(E^u)}\subset\Delta\,,$$
 obtaining the desired result.

For this, is enough to note that  $ \Lambda $  is in the hypothesis of the lemma \ref{lema3}.
Therefore as $\pi^s(\Gamma_{\delta})$  is  dense  in  $E^u$ for a
$ \delta $ sufficiently small, we can join  $ x_0 $ with itself by a curve $ \delta $-adapted such
that when lifted,  it links $x_0 $
with any point of $ \Gamma_{\delta}$. For an appropriate$ \delta $ , and any $ x \in \Gamma_{\delta}$, we have that $x_0$ and $x$
 are n-$ \delta $-related for
some $ n $, as desired.
\endproof


\section{Mañe's robustly transitive diffeomorphisms that is not Anosov}

In this section we will describe  an example constructed  by Mañe in \cite{M}. This example is
very well described in numerous references (see for instance \cite{BDV},  or \cite{PS}),
but we will include a description for the convenience of the reader, and because we will emphasize some
properties of the example that will be useful later on. However we will not include the proofs,which can be
found in any of the given references.

As in the previous section, let us starts with a linear
 Anosov diffeomorphism $f_A$ in $\mathbb{T}^3$ that is
induced form  $ A \in GL(3, \Z)$ which is a hyperbolic toral automorphism with only one
eigenvalue grater than one, and
all eigenvalues real, positive, simple, and irrational.
Let $0<\lambda^s< \lambda^c <1<\lambda^u$ be the eigenvalues.
Let $\mathcal{F}^c$
be the foliation corresponding to the eigenvalue $\lambda^c$, similarly with $\mathcal{F}^s$
and $\mathcal{F}^u$. We remind you that all of these leaves are dense.
 We may also assume that $f_A$ has at least two
fixed points, $x_0$ and $x_1$, and that unstable eigenvalue $\lambda^u$, have modulus greater than 3
(if not, replace $A$ by some power).

Following the construction in \cite{M} we define $f$ by modifying $f_A$ in a sufficiently small domain $C$
contained in $B_{\frac{\rho}{2}}(x_1)$ keeping invariant the foliation $\mathcal{F}^c$. Where $\rho> 0$ is
a small enough number to be determined in what follows.
Let us observe that $f_A|_{B_{\frac{\rho}{2}}(x_1)^c}=f|_{B_{\frac{\rho}{2}}(x_1)^c}$. In particular
\begin{equation}\label{invariantes}
\Gamma=\bigcap_{n\in\Z}f^n(B_{\frac{\rho}{2}}(x_1)^c)=\bigcap_{n\in\Z}f_A^n(B_{\frac{\rho}{2}}(x_1)^c)\,.
\end{equation}
\begin{figure}
\begin{center}
\scalebox{0.35}{\includegraphics{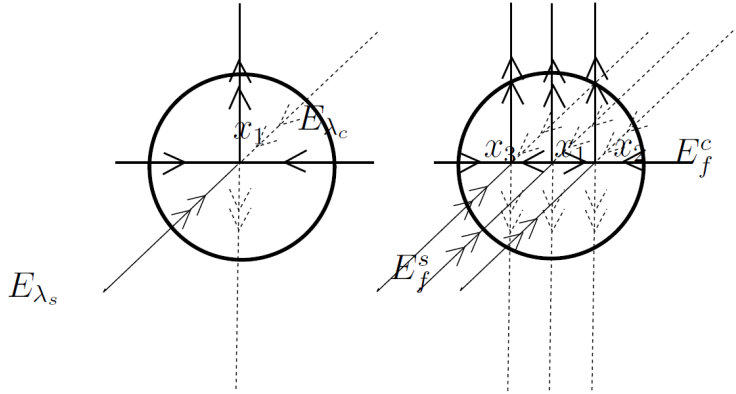}}
\caption{Perturbing a neighborhood of $x_1$}\label{pert}
\end{center}
\end{figure}

 We can take $\rho$ sufficiently small so that $x_0\in\Gamma$.
 Inside $C$ the point $x_1$ undergoes a bifurcation as shown in the figure \ref{pert}, in the direction of $\mathcal{F}^c$,  which changes the unstable index
 of $x_1$ increasing it in 1. Also two other fixed points, $x_2$ and $x_3$ are created, with the same index $x_1$  had under
 $f_A$.

 As a result,  we get a difeomorphism $f$ which is strongly partially hyperbolic. That is
 $$ T\mathbb{T}^3 = E_f^s\oplus E_f^c\oplus E_f^u\,,$$
where $E_f^s$ is uniformly contracting and $ E_f^u$ is uniformly expanding.
In fact, $E_f^s$ and $E_f^u$ are contained in some small cones around $E^s$ and  $E^u$
respectively. Then by a well known results (see \cite{hps}) we get that the bundles
 $E_f^s$ and $E_f^u$ are uniquely integrable  to foliations $F_f^s$ and
$F_f^u$  called the (strong)
stable and unstable foliations.
Moreover, they are quasi-isometric
Since we preserved the central foliation we have   $\mathcal{F}_f^c=\mathcal{F}^c$,
 $E_f^s\oplus E_f^c$ and $E_f^u\oplus E_f^c$ are also uniquely integrable, by what we call the center-stable and center-unstable
foliations respectively.
In \cite{M} it is shown that the leaves of $\mathcal{F}_f^c$ are dense in $ \mathbb{T}^3$ (see also \cite{BDV}),
and also in a robust fashion.

It is particularly relevant for us that, not onely is the central foliation minimal, but also the
unstable foliation is minimal as well.
This is shown  for instance in the following theorem from \cite{PS} (page 5).
\begin{theorem}\label{minimalu} (2.0.1 in \cite{PS}).  There exists a   neighborhood $\mathcal{U}$ of $f$, in the $C^1$ topology such that
for every $g\in\mathcal{U}$ the bundles $E^c_g$, $E^s_g$ and $E^u_g$,   uniquely integrates
 to  invariant foliations ($\mathcal{F}^c_g$, $\mathcal{F}^s_g$ and $\mathcal{F}^u_g$, respectively).
 Furthermore, the central and unstable
foliations $g \in\mathcal{U} $ are minimal, i.e., all leaves are dense.
\end{theorem}


The following lemma is a consequence of the shadowing theorem (see \cite{S}).
\begin{lemma}\label{H}
Let $ A \in GL(3, \Z)$ which is a hyperbolic toral automorphism and let $G : \R^3\to\R^3$
be a homeomorphism such that $\norm{A(x)-G(x)} \leq r$ for all $x \in \R^3$. Then there exists
$H : \R^3\to\R^3$
continuous and onto such that $A \circ H = H \circ G$. Moreover
$\norm{H(x) - x} <C.r$  for al $x$.
\end{lemma}


Note that $H(x) = H(y)$ if and only if $\norm{G^n(x) - G^n(y)} \leq 2Cr\,\,\,\,\forall n \in \Z$.
This is a consequence of the uniqueness in the shadowing theorem.

Since
$G$ is isotopic to $A$, 
$H$ induces an  $h: \mathbb{T}^3 \to \mathbb{T}^3$
continuous
and onto such that $f_A \circ h = h \circ g $ and $d_{C_0}(h, id) = rd_{C_0}(f_A, g)$.

As a consecuence of this we have:
\begin{lemma}\label{Huc}
With the above notation, $H : \R^3\to\R^3$ is uniformly continuous  for every $x$.
\end{lemma}



Now let us see how $H$ behaves with respect of the invariant foliations.
\begin{lemma}\label{Hprop}
For $H$, $A$ and $G$ as above we have that
\begin{enumerate}
\item  $H\widehat{\mathcal{F}}^{cu}_G(x)=\widehat{\mathcal{F}}^{cu}_A(H(x))$ and
$H\widehat{\mathcal{F}}^{cs}_G(x)=\widehat{\mathcal{F}}^{cs}_A(H(x))$.
\item $H\widehat{\mathcal{F}}^{c}_G(x)=\widehat{\mathcal{F}}^{c}_A(H(x))$.
\item $H\widehat{\mathcal{F}}^{u}_G(x)=\widehat{\mathcal{F}}^{u}_A(H(x)) = H(x)+E_A^u$
and $H\mid_{\widehat{\mathcal{F}}^{c}_G(x)}$ is a homeomorphism for every $x$.
\item For any $x, y \in \R^3$,
$$ \#\set{\widehat{\mathcal{F}}^{cs}_G(x) \cap \widehat{\mathcal{F}}^{u}_G(y)=1 }\text{ and }
\#\set{\widehat{\mathcal{F}}^{cu}_G(x) \cap \widehat{\mathcal{F}}^{s}_G(y) }= 1\,.$$
\item If $ H(x) = H(y)$, Then $x$ and $y$ belong to
the same central leaf.
\end{enumerate}
\end{lemma}
These results follow mainly from the expansivity of $A $ and the fact that $\norm{H(x) - x} < Cr$. For
 a proof see \cite{PS}.
It can also be shown that $h:\mathbb{T}^3 \to \mathbb{T}^3 $ inherits similar properties.



\section{Proof of Theorem \ref{main2}: a set that is robustly not premaximal in $\mathbb{T}^3$}



Let $f:\mathbb{T}^3\to\mathbb{T}^3$ be as in the previous section, the diffeomorphism form Mañe's example, and let
us consider a $C^1$ ball around $f$, $\mathcal{U}$.
In this section we will prove that for any $g\in \mathcal{U}$,  there is a set on
 $\mathbb{T}^3$ that cannot be included on any set with local product structure.

For this we will show that the set $\Lambda$ from section \ref{sec2} does not intersect some ball around $x_1$
So possibly  taking a smaller $\rho$ we can construct a diffeomorphism $f$ as the one from the previous
section and such that
$\Lambda\subset\Gamma=\bigcap_{n\in\Z}f^n(B_{\frac{\rho}{2}}(x_1)^c)\,.$

Note that the set $\Gamma=\bigcap_{n\in\Z}f^n(B_{\frac{\rho}{2}}(x_1)^c)\,$ can be made to be transitive.

So $\Lambda$ is a compact, hyperbolic set, invariant under $f$ since it is invariant under $f_A$ and by equation (\ref{invariantes}).
For any $g $ sufficiently close to $f$, there is a hyperbolic set $\Lambda_g$ which is the hyperbolic continuation
of $\Lambda$, and that has essentially the same properties in all that concerns us. We will call both sets
 $\Lambda$, for simplicity.
We aim to prove that if there is a set $\Delta $ containing  $\Lambda$ with local product structure, then
$\Delta\cap\mathcal{F}_g^u(x_0)$ is dense in some small interval of $\mathcal{F}_g^u(x_0)$, and then  $\Delta $
is dense in $\mathbb{T}^3$ in virtue of the minimality of $\mathcal{F}_g^u$ (\ref{minimalu}). This is a contradiction since $g$
is not Anosov.
Since in this context the unstable leaves are not parallel it would be convenient to redefine the
n-$ \varepsilon $-relation.

Let $ p^{cs}_x: \R^3 \to \widehat{\mathcal{F}}_{G}^{u}(x) $ and $ p^{u}_x: \R^3 \to\widehat{\mathcal{F}}_{G}^{cs}(x)$
be the projections along the  center stable and  unstable foliation respectively.
We note as  $ \widehat{\Delta} $  the lift of $ \Delta $.

\begin{definition}
We say that $ x $ and $ y $ are n-$ \varepsilon $-related in $\widehat{\Delta}$ if there exists
 sequence of point $ x = x_0, x_1, \dots, x_n = y $ such that:

\begin{itemize}
\item $ x_i \in \widehat{\Delta} $ for $ i = 1, \dots, n $
\item $ d(p^{u}_{x_ {i}}(x_ {i +1}),x_ {i+1}) \leq \varepsilon \text{ for } 1 \leq i \leq n $
\item $  d(p^{cs}_{x_ {i+1}}(x_ {i }),x_ {i}) \leq \varepsilon $ for $ 1  \leq i \leq n $
\end{itemize}
\end{definition}

The main problem  which we are dealing with now, is that the lemma (\ref{EPG})
relies heavily on the linearity of $ A $. We will fix this problem by  finding a tube $V$ around $(0,0,0)$ so that both the distance between
the center-sable foliations of $x$ and $y$
and the distance between the unstable foliations in $V$ are small  when $x$ and $y$ are close enough.
The interval of the unstable foliation  in which
$ \widehat{\Delta}\cap\widehat{\mathcal{F}}_G^u((0,0,0))$  will dense, will be contained in this $V$.

Another important difference is that \ref{remark} also makes a strong use of the linearity therefore we will not try to prove that the projection of all $\Gamma_{\delta}$ is in $\widehat{\Delta}$. It will be enough to find a point of $\Gamma_{\delta}$ outside $V$ and project the points of the $\delta$-chain joining $(0,0,0)$ with that point.

For two points $x$ and $y$ in the same leaf of the unstable foliation, we define $l^u(x,y)$ to be the
 length of the arc joining $x$ with $y$.
For a fixed $\varepsilon$, we will prove first that for any tow points $x$, $y$ in $\R^3$, there exist  a $\delta$ such that if $d(x,y)< \delta$.
Then, if we choose any $z$ in $\widehat{\mathcal{F}}^{cs}(x)$, then $l^u(z,p^u_y(z))<\varepsilon$.

\begin{lemma}\label{arreglandoV}
For any $\varepsilon>0$ there exists
a $\delta$ such that for every $x$ and   $y\in \widehat{\mathcal{F}}^{u}_G(x)$ such that $l^u(x,y)\leq\delta$ then   $l^u(z,p^u_y(z))<\varepsilon$,
for any $z$ in $\widehat{\mathcal{F}}_{G}^{cs}(x)$.
\end{lemma}
\proof
Suppose that this is not the case. Then there must exist an $\varepsilon_0$ such that there exist there
sequences $\set{x_n}_{n\in\N}$, $\set{y_n}_{n\in\N}\subset\widehat{\mathcal{F}}^{u}_G(x)$ and $\set{z_n}_{n\in\N}$ such that $l^u(x_n,y_n)\leq 1/n$,
 and $l^u(z_n,p^u_y(z_n))\geq\varepsilon_0$.

 \begin{figure}[h]
\begin{center}
\scalebox{0.30}{\includegraphics{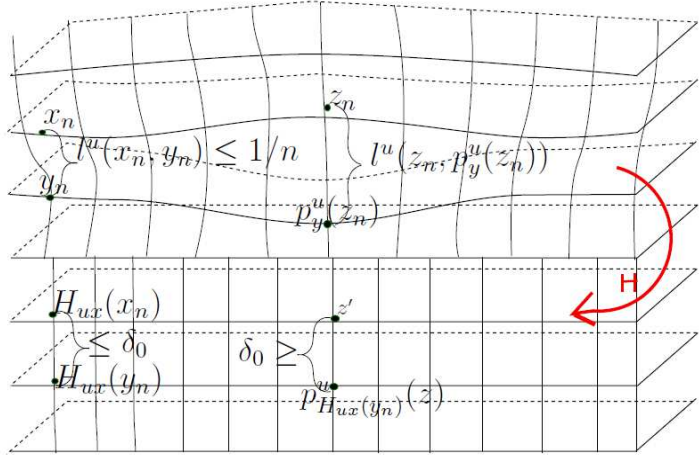}}
\caption{$H$ acting on the foliations }\label{hola}
\end{center}
\end{figure}
Now let us recall that from (\ref{Hprop}) we have that $H\mid_{\widehat{\mathcal{F}}^{u}_G(x)}$ is a
homeomorphism, for simplicity we note $H\mid_{\widehat{\mathcal{F}}^{u}_G(x)}=H_{ux}$.
Let $\delta_0$ be the one given by de uniform continuity of $H$ (\ref{Huc}).  Since $H_{ux}$
is a homeomorphism,  for $\delta'$ we can find a $\delta_0$ (independent of $x$) such that if $x$ and $y$ are such that  $y\in \widehat{\mathcal{F}}^{u}_G(x)$ and $d(H_{ux}(x),H_{ux}(y))\leq\delta_0$, then
$l^u(x,y)\leq\delta'$.

Let us consider $n_0$ such that $1/n_0<\delta'$ and  $z'= H(z_{n_0})$. Note
that $$z'\in\widehat{\mathcal{F}}_{A}^{cs}(H(x_{n_0}))\,.$$

For perhaps a bigger $n$, we  have that   $l^u(x_n,y_n)\leq 1/n$, and $d(H(x_n),H(y_n))\leq\delta_0$, from the continuity of $H$.
But for $A$, $\widehat{\mathcal{F}}_{A}^{cs}$ are parallel planes so since
the length of the unstable segment between
$z'$ and $p'^u_{H(y_{n})}(z')$ is less than $\delta_0$ (see figure \ref{hola}) and therefore
$$\varepsilon_0> \delta'>l^u(z_n,H_{ux}^{-1}(p'^u_{H(y_{n})}(z')))=l^u(z_n,p^u_y(z_n))\geq\varepsilon_0\,.$$
\endproof
 If two points are sufficiently close their unstable manifolds remain close in some neighborhood. This is a consequence of the continuity of the foliation.
 \begin{figure}[h]
\begin{center}
\scalebox{0.35}{\includegraphics{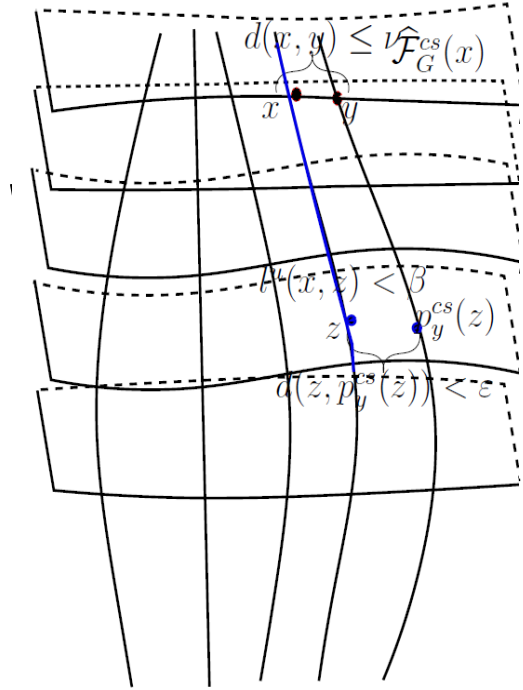}}
\caption{Separation between the unstable leafs  }\label{sep}
\end{center}
\end{figure}

For every  $\varepsilon>0$,
there is  $\beta>0 $ and  $\eta>0 $
such that if $y\in\widehat{\mathcal{F}}_{G}^{cs}(x)$
and $d(x,y )<\eta$, then for any $z\in\widehat{\mathcal{F}}_{G}^{u}(x)$ such that $l^{u}(x,z )<\beta$,
we have that $d(z,p^{cs}_{y}(z) )<\varepsilon$.
We can also take $\beta$ to be uniform since the foliations
are lifts of foliations in a compact set (see figure \ref{sep}).

Now  we will put everything together.
Let $\varepsilon=\delta_p$ be the one from the local product structure of $\widetilde{\Delta}$. For this $\varepsilon$ we find $\eta>0$
and $\beta$ from our previous observation. This will ensure us that if $d(x,y )<\eta$ their unstable leaves will remain closer than  $\delta_p$  in a ball of radius $\beta$ from $x$.

We can choose the $\varepsilon_0$ from the lemma (\ref{arreglandoV})
 smaller $\delta_p$ and $\beta$,
 So the lemma ensures us that there exists a
 $\delta_0$ such that if  $x$ and   $y\in \widehat{\mathcal{F}}^{u}_G(x)$ and $l^u(x,y)\leq\delta_0$ then  the center stable foliations of $x$ and $y$ will not separate more than $\varepsilon_0$.

  For this last $\delta_0$ we will take a compact neighborhood of $(0,0,0)$ in $\widehat{\mathcal{F}}^{u}((0,0,0))$ that we call $U^u$ with $diam(U^u)<\delta_0$.
In this conditions we define $$V=\bigcup_{x\in U^u}\widehat{\mathcal{F}}_{G}^{cs}(x)\,.$$
We have proved the following for $V$.

 \begin{figure}[h]
\begin{center}
 \scalebox{0.35}{\includegraphics{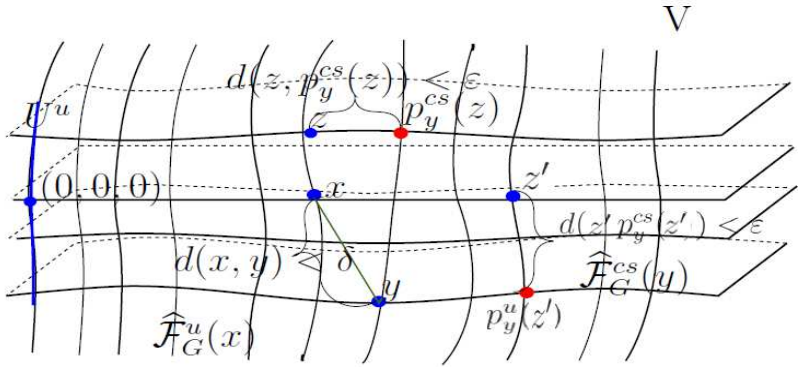}}
\caption{The tube $V$ }\label{V}
\end{center}
\end{figure}

\begin{lemma}\label{defV}

There exist a compact neighborhood of $(0,0,0)$, $U^u$, in $\widehat{\mathcal{F}}^{u}((0,0,0))$, such that the set $V$
define as $V=\bigcup_{x\in U^u}\widehat{\mathcal{F}}_{G}^{cs}(x)\,,$ satisfies:
\begin{itemize}
\item  If $x\in V$ and $ y\in\widehat{\mathcal{F}}_{G}^{u}(x)\cap V$, then $d(x,y)<\delta_p $  and $d(x,y)<\beta$.
\item For every $\varepsilon$  there is $\delta$ such that if $x\in V$ and $d(x,y)<\delta$ then
$d(z,p^u_y(z))<\varepsilon$,
for any $z$ in $\widehat{\mathcal{F}}_{G}^{cs}(x)$.
\item   For every $\varepsilon$  there is $\delta$ such that if $x\in V$ and $d(x,y)<\delta$ then
$d(z,p^{cs}_y(z))<\varepsilon$,
for any $z$ in $\widehat{\mathcal{F}}_{G}^{u}(x)\cap V$.
\end{itemize}
\end{lemma}

The following implies Theorem \ref{main2}:
\begin{theorem}
Let $f:\mathbb{T}^3\to \mathbb{T}^3$ be as in Mañe's example and let $g\in \mathcal{U}(f)$, be a
 difeomorphism sufficiently close to $f$. There exists
 $\Lambda$
a compact, connected,non trivial, $g$-invariant, hyperbolic set such that  $x_0\in\Lambda$
and $\Lambda$ is not included in any set with local product structure.
 \end{theorem}
\proof
 We first recall that for every   $g\in \mathcal{U}(f)$ the  existence of a set $\Lambda$ which is a compact, connected, non trivial, $g$-invariant, hyperbolic set, has already been stated.

As before we start by supposing that there is a set $\Delta$ with local product structure
containing $\Lambda$.
Recall that the strategy to see that such a $\Delta$  can not exist,  is to find an interval of the unstable foliation  in which
$ \widehat{\Delta}\cap\widehat{\mathcal{F}}_G^u(x_0)$  will dense.

Let $V\subset \R^3$ be an open tube as defined in \ref{defV},

\begin{af}\label{epg2}
Suppose that $x_0,\dots,x_n\,\in \widehat{\Delta}\cap V$ are n-$\varepsilon$-related, and :
\begin{itemize}

\item For any $z$ in $\widehat{\mathcal{F}}_{G}^{cs}(x_{i+1})$  ,
$d(z,p^u_{x_{i}}(z))<\varepsilon$, for every $0\leq i\leq n-1$.
\item For any  $z$ in $\widehat{\mathcal{F}}_{G}^{u}(x_i)\cap V$,
$d(z,p^{cs}_{x_{i+1}}(z))<\varepsilon$,for every $0\leq i\leq n-1$.
\end{itemize}
Then $p^{cs}_{x_{0}}(x_i)\in\widehat{\Delta}$ for every $0\leq i\leq n-1$.
\end{af}

\proof
We will prove this by induction. The base case is given by the fact that we chose  $0 <\varepsilon$ smaller
than $\delta_p/2$, where $\delta_p$ is from the local product structure  of $\widehat{\Delta}$.
Suppose now that $ x, y \in \widehat{\Delta}\cap V $ are
n-$\varepsilon$-related. We have $ x = x_0, x_1, \dots, x_n = y $ as in the definition. We define
$ \overline{x_j} = \widehat{\mathcal{F}}_G^{u}(x_j) \cap \widehat{\mathcal{F}}_G^{cs}(x_ {j +1}) $ for $ 0 \leq j \leq n-1 $. Note that
$ \overline{ x_0} $ and $ \overline{x_ {n-1}} $ are (n-1)-$\varepsilon$-related because:
\begin{itemize}
\item $\overline{ x_j} \in \widehat{\Delta} $ for $ j = 1, \dots, n-1 $ by induction hypothesis,
\item Since $\overline{x_j} $ is in $\widehat{\mathcal{F}}_{G}^{cs}(x_{j+1})$, by our hypothesis we have that
$d(\overline{x_j},p^u_{x_{j}}(\overline{x_j}))<\varepsilon$
On the other hand since $ \overline{x_{j-1}} \in \widehat{\mathcal{F}}_G^{cs}(x_ {j }) $ then,
$d(\overline{x_j},p^u_{\overline{x_{j-1}}}(\overline{x_j}))<\varepsilon$.
\item  Since $ \overline{x_j} = \widehat{\mathcal{F}}_G^{u}(x_j) \cap \widehat{\mathcal{F}}_G^{cs}(x_ {j +1}) $
, we have that, $ \overline{x_j}\in\widehat{\mathcal{F}}_{G}^{u}(x_j)\cap V$. Again by our assertion´s hypothesis we have that
$d(\overline{x_j},p^{cs}_{x_{i+1}}(\overline{x_j}))<\varepsilon$
On the other hand,  $ \overline{x_{j+1}} = \widehat{\mathcal{F}}_G^{u}(x_ {j +1}) $, so
$d(\overline{x_j},p^{cs}_{\overline{x_{j+1}}}(\overline{x_j}))<\varepsilon$.
\end{itemize}

   This allows us to conclude that  $p^{cs}_{x_{0}}(x_i)\in\widehat{\Delta}$.
This proves our assertion.

\endproof

Returning to the proof of the theorem,
for any $\varepsilon<\delta_p/2$ we take the $\delta<\varepsilon$ from the definition of $V$ as in lemma (\ref{defV}).

Let us take any point $q$ of $\Gamma_{\delta}$ which is not in $V$.
 $(0,0,0)$ is  n-$\delta/2$-related to $q$. We call $x_{j+1}$ the firs element of the sequence relating $(0,0,0)$
 to $q$ that is not in $V$. As in the proof of the assertion, we can construct a new sequence from the
  sequence that  $j$-$\varepsilon$-relates $(0,0,0)$ to $x_j$ as follows. 

   We define
$ \overline{x_i} = \mathcal{F}^{u}(x_i) \cap \mathcal{F}^{cs}(x_ {i +1}) $ for $ 0 \leq i \leq j-1 $. Note that
$ \overline{ x_0} $ and $ \overline{x_ {j-1}} $ are $(j-1)$-$\varepsilon$-related because:
\begin{itemize}
\item $\overline{ x_i} \in \widehat{\Delta} $ for $ i = 1, \dots, j-1 $ since
$\delta/2<\varepsilon<\delta_p/2$ from the local product structure.
\item   $ \overline{x_i}\in \mathcal{F}^{cs}(x_ {i +1}) $ , and $x_i,x_{i+1}\in V$ with $d(x_i,x_{i+1})<\delta$.
From \ref{defV} we have that
$d(\overline{x_i},p^u_{x_i}(\overline{x_i}))<\varepsilon$. As in the assertion this implies
$d(\overline{x_i},p^u_{\overline{x_{i-1}}}(\overline{x_i}))<\varepsilon$ for $ i = 1, \dots, j-1 $.

\item Since $ \overline{x_i} = \widehat{\mathcal{F}}_G^{u}(x_i) \cap \widehat{\mathcal{F}}_G^{cs}(x_ {ji+1}) $
, we have that, $ \overline{x_i}\in\widehat{\mathcal{F}}_{G}^{u}(x_i)\cap V$. Since $d(x_i,x_{i+1})<\delta$. As before from \ref{defV} we have that

$d(\overline{x_i},p^{cs}_{x_{i+1}}(\overline{x_j}))<\varepsilon$
On the other hand,  $ \overline{x_{j+1}} = \widehat{\mathcal{F}}_G^{u}(x_ {i +1}) $, so
$d(\overline{x_i},p^{cs}_{\overline{x_{i+1}}}(\overline{x_j}))<\varepsilon$ for $ i = 1, \dots, j-1 $.
\end{itemize}
\begin{figure}
\begin{center}
\scalebox{0.35}{\includegraphics{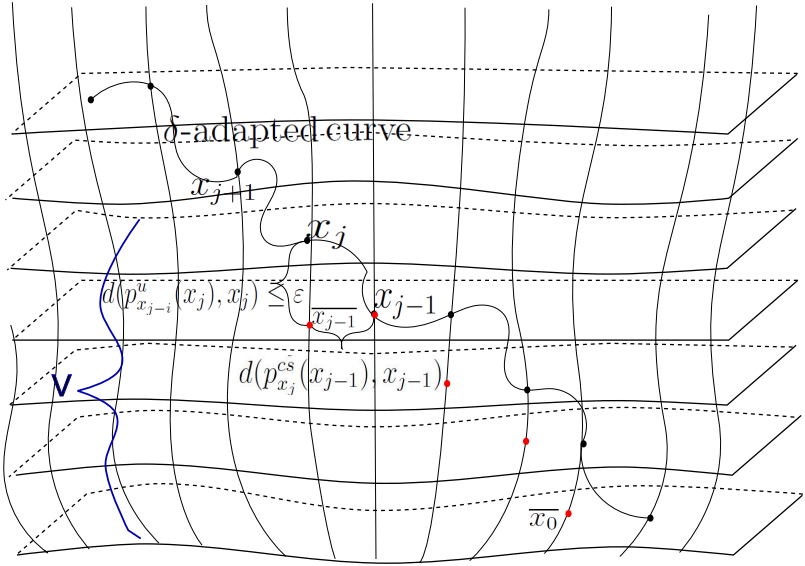}}
\caption{A n-$ \varepsilon $-relation between  $ (0,0,0) $ and $ x_j $.  }\label{tda}
\end{center}
\end{figure}
In addition to this, $\widehat{\mathcal{F}}^{u}(\overline{x_i}) \cap V=\widehat{\mathcal{F}}^{u}(x_i)
 \cap V$ for$ 1 \leq i \leq j-1 $ and $\widehat{\mathcal{F}}_{G}^{cs}(\overline{x_i})=\widehat{\mathcal{F}}_{G}^{cs}(x_{i+1})$ . Recall that since  $x_i\in V$ for $ i = 1, \dots, j$ and
 $d(x_i,x_{i+1})\leq\delta$,
 so by \ref{defV} we have:
 \begin{itemize}

\item For any $z$ in $\widehat{\mathcal{F}}_{G}^{cs}(\overline{x_i})$  ,
$d(z,p^u_{x_{i}}(z))=d(z,p^u_{\overline{x_{i-1}}}(z))<\varepsilon$, for every $0\leq i\leq j-1$.
\item For any  $z$ in $\widehat{\mathcal{F}}^{u}(\overline{x_i})\cap V$,
$d(z,p^{cs}_{x_{i+1}}(z))=d(z,p^{cs}_{\overline{x_{i+1}}}(z))<\varepsilon$,for every $0\leq i\leq j-1$.

\end{itemize}
 This proves that the new sequence we have defined is in the hypothesis of our assertion so,
  $p^{cs}_{\overline{x_{0}}}(\overline{x_i})\in\widehat{\Delta}$ for $i=0,\dots,j-1$.
  But then $p^{cs}_{x_{0}}(x_i)\in\widehat{\Delta}$ for $i=0,\dots,j$.

Let $ U^{u+} $ and $ U^{u-} $, be the positive  and negative sub intervals of $ U^{u} $.
Arguing by contradiction we suppose that $ \widehat {\Delta} \cap U^{u} $ is not dense in any of these sub intervals.
 If this is so, there must be some gaps of size at least $\gamma$ in each sub interval,  for which
 there are no points of $ \widehat {\Delta} \cap U^{u} $ in these gaps.
 But choosing $\varepsilon\leq\gamma/2$, $\delta$ for this $\varepsilon$ and $q$  a point in the corresponding $\Gamma_{\delta}$, $(0,0,0)$ is n-$\delta$-related to $q$.

  Therefore for $j $ as before, $p^{cs}_{x_0}(x_i)\in\widehat{\Delta}$ for $i=0,\dots,j$.
  Since $x_{j+1}$ is out of $V$, one of the subintervals, $ U^{u+} $ or $ U^{u-} $ has at least
   one point of $\widehat{\Delta}$ in every gap of size $2\varepsilon$ leading to a contradiction with
    our previous assumption.
   Iterating $U^u$ and since $\widehat {\Delta}$ is invariant, we have then that
    $$\bigcup_{n\in\Z^3}g^n(U^u)=E^u\subset\widehat {\Delta}\text{ and }\mathbb{T}^3 =\overline{\pi(E^u)}\subset\Delta\,.$$
    \endproof

  \section{in dimension 2  transitive sets are locally premaximal}

In this section we prove Theorem \ref{main3}, allowing us to complete the answer to the  question of whether or not it is possible to construct transitive sets which are not locally premaximal in dimensions smaller than 4.
 For this, we rely on a result by Anosov \cite{A1} (for a proof in an other context see also \cite{BG} proposition 4.3) and one by Fisher in \cite{Fi}, which we state below.
 \begin{theorem}(Anosov) \label{Anosov}
 Let $ f: M \to M $ be a diffeomorphism of a compact manifold and let $ F \subset M $ be a compact hyperbolic invariant set with zero topological dimension.
 For every neighborhood $ U $ of $ F $, there exists a set $ \Lambda $ which is compact invariant and hyperbolic and has local product structure, such that $ F \subset \Lambda \subset U $.
 \end{theorem}

 \vspace {-0.6cm}

 \begin{theorem} (Fisher) \label{fisher3}
Let $ f: M \to M $ be a diffeomorphism of a compact surface $ M $. If $ \Lambda $ has local product structure and has non empty interior, then $ \Lambda = M = \mathbb{T}^2 $.
 \end{theorem}


\begin{Remark}
Since $M$ is a Hausdorff locally compact space, (we will have a compact surface actually), then zero dimensional subsets are exactly the totally disconnected subsets.
\end{Remark}
We will prove some lemmas that will imply Theorem \ref{main3}.

In what follows $ f: M \to M $ will be  a diffeomorphism, of $ M $ a compact surface  and $ \Lambda \subset M $ will be  a compact hyperbolic invariant set such that $ \Omega (f \mid_{\Lambda}) = \Lambda $

  \vspace {0.2cm}

\begin{definition}
Let $ \Lambda \subset M $  be  a compact hyperbolic invariant set. We note $ \Lambda_0$ is the union of all points $p$ in $\Lambda$ such that the connected component of $p$ in $\Lambda$ is $p$ itself. We define $\Lambda_1 $ as $\Lambda_1 =\Lambda\setminus\Lambda_0$.
\end{definition}

 Note that $ \Lambda_0$  is totally disconnected and therefore $0$ dimensional.


As before, for a sufficiently small neighborhood of $ x \in \Lambda_1 $ we define $ p^s_ {x} :U (x) \cap \Lambda_1\to W^u_{loc}(x)  $ to be the local projection along the stable manifolds. We define analogously  $ p^u_ {x} :U (x)\cap \Lambda_1\to W^s_{loc}(x)  $. Recall that both stable and unstable manifolds are one dimensional.

  \vspace {0.2cm}

\begin{lemma}\label{periodens}
 Periodic points are dense in $ \Lambda_1 $. Moreover for any $x\in \Lambda_1 $ we have that $W^u_{loc}(x)\subset \Lambda_1$ or $W^s_{loc}(x)\subset \Lambda_1$ (or both).
  \end{lemma}
  \proof
We note the local connected component of $x$ in $ \Lambda_1$ as $lcc(x)$ and from the definition of $ \Lambda_1$ we have that  $lcc(x)$ is not trivial. Then for any  $ x \in \Lambda_1 $ we have that either $ p^s_ {x} (lcc(x))$ contains a nontrivial connected set (an arc) or
$ p^u_ {x} (lcc(x))$  contains a nontrivial connected set.
This will also be true for any smaller $U (x)$.

 Since $ \Omega (f \mid_{\Lambda}) = \Lambda $ , using the shadowing theorem we have that if $ x \in \Lambda $ then it is approximated by periodic points (which a priori would not be in $ \Lambda $ ).

 Let us suppose that $ p^s_{x} (lcc(x)) $ contains an arc. Then  we define $$  V =\bigcup_{y \in lcc(x)} W^s_{loc} (y) $$
and we have that $\mathring{V}$ is open and not empty. We take $z\in\mathring{V}\cap lcc(x)$ and  $p^s(z)$ is in the interior of $ p^ s_ {x} (lcc(x)) $ in the relative topology of $W^u(x)$. If $ \set{p_n}_{n \in \N} $ is such that $ p_n \to z$   (where $p_n$ are periodic points), then $ p_n \in \mathring{V} $ for all $ n $ greater than some $ n_0 $.

\begin{figure}[h]
\begin{center}
\scalebox{0.35}{\includegraphics{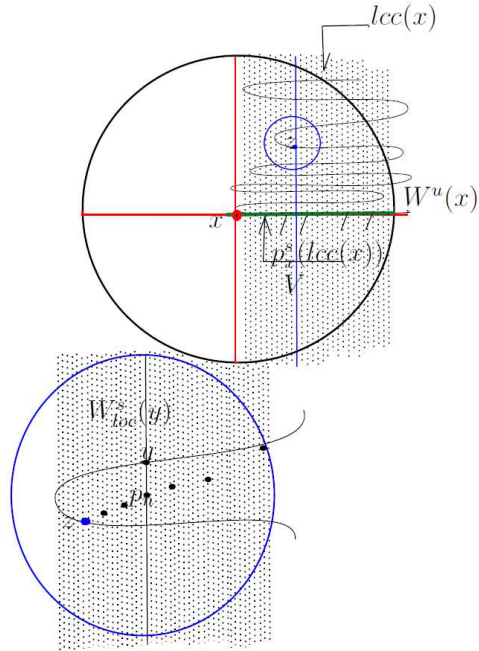}}
\caption{ the local connected component of $x$ and $ p^s_{x} (lcc(x)) $ }
\end{center}
\end{figure}

Since $ \Lambda $ is invariant, and $ f $ a diffeomorphism it follows that $ \Lambda_1$ is invariant too ( otherwise  it would take a nontrivial connected component into a point). This implies that $ \omega (p_n) \in \Lambda_1 $ but then $ p_n \in \Lambda_1 $ for all $ n> n_0$.

 \vspace {0.2cm}
Now we take $ p_n $ sufficiently close to $ z$ such that $p^s_x(p_n)$ is in the interior of $ p^ s_ {x} (lcc(x)) $ in the relative topology of $W^u(x)$. Then $ W^s_{loc} (p_n) \cap lcc(x) \neq \emptyset $.

 Iterating for the future $ f^n (lcc(x)) $ accumulates on the unstable manifold containing $ p_n $ in its interior  and since $ diam (lcc(x))> 0 $, $ f^n (lcc(x)) $ accumulates on an arc of $ W_{loc}^u (p_n) $. Therefore this arc of  $ W_{loc}^u  (p_n) $ is contained in $ \Lambda $, since $ \Lambda $ is compact. Such an arc must be  contained in $ \Lambda_1 $ because it clearly does not belong to $\Lambda_0 $ . The invariance of $ \Lambda_1 $ implies that $ W^u (p_n) \subset \Lambda_1 $ and as $ p_n \to z $ implies that $ W^u (z) \subset \Lambda_1 $. We can take now a sequence of $z_n$ as before that are contained in neighborhoods $ U_n (x)$ which are each time smaller, and $z_n\to x$ so $ W^u (x) \subset \Lambda_1 $.


The situation is analogous if $ p^u_{x} (lcc(x)) $ is  a stable arc.
  \endproof

  \vspace {0.3cm}

  \begin{lemma}\label{masconj}
 Let $ \Lambda_1^u=\set{x\in\Lambda_1 \tq W^u (x) \subset \Lambda_1 }$ and $ \Lambda_1^s=\set{x\in\Lambda_1 \tq W^s (x) \subset \Lambda_1 }$. Then $\Lambda_1^u$ and $\Lambda_1^s$ are compact sets with local product structure.
  \end{lemma}

\proof
 We will prove first that $\Lambda_1^u$ is closed and therefore compact. Let $\set{x_n}_{n\in\N}\subset \Lambda_1^u$   and $ x_n \to y$   then the unstable manifold  of $ y $  is in $ \Lambda $ since $\Lambda$ is compact. Hence $ y \in \Lambda_1^u $.
 Let $ x, y \in \Lambda_1^u$ such that $ d (x, y) <\delta $ for some $ \delta $ appropriate  such that $ W^u_ {loc} (x) \cap W^s_{loc} (y)\neq \emptyset$. Since $ W^ u (x) \subset \Lambda_1 $, and $ W^u (y) \subset \Lambda_1 $, $$ W^u_ {loc} (x) \cap W^s_{loc} (y) \in W^u (x) \subset \Lambda_1 \,. $$

The situation is analogous for $\Lambda_1^s$ .
\endproof

  \vspace {0.2cm}

\begin{corollary}\label{cor1}
The set $\Lambda_1$ is compact, invariant and has local product structure.
\end{corollary}
\proof
The sets  $ \Lambda_1^u$ and $ \Lambda_1^s$ defined in lemma \ref{masconj} are such that
 $$ \Lambda_1=\Lambda_1^u \cup \Lambda_1^s\,,$$ as a consequence of lemma \ref{periodens}. Therefore  $\Lambda_1$ is compact and has local product structure.
\endproof

  \vspace {0.2cm}

\begin{corollary}\label{forma}
Either $\Lambda_1$  is the disjoint union of an attractor $\Lambda_1^s$ and a repeller $\Lambda_1^u$, or $\Lambda=M=\mathbb{T}^2$ and $f$ is Anosov.
\end{corollary}

\proof
Suppose that there is $x\in\Lambda_1^u \cap \Lambda_1^s$. Then $ W^u (x)$ and $ W^s(x) \in  \Lambda_1  $. Since $\Lambda_1  $ has local product structure this implies that $$\bigcup_{y\in W^s_{loc} (x)} W^u_{loc} (y)\subset \Lambda_1 \,, $$ and so $\Lambda_1  $ has non empty interior. It follows from Theorem \ref{fisher3} that $\Lambda=\Lambda_1=M=\mathbb{T}^2$ and $f$ is Anosov.
\endproof

\vspace {0.2cm}

\begin{lemma}\label{ultimo}
 The set $\Lambda_0$ is compact and disjoint from $\Lambda_1$. Moreover for any neighborhood $V$ of $\Lambda_0$, we can find a set $\Lambda_0'$ with local product structure such that $\Lambda_0'\cap\Lambda_1=\emptyset$ and $ \Lambda_0 \subset \Lambda'_0 \subset V $.
  \end{lemma}

\proof
From its definition, $\Lambda_0$ is  disjoint from $\Lambda_1$. Suppose that $\Lambda_0$ is not empty (the lemma holds trivially if it is empty).
Suppose that there is a sequence $\set{x_n}_{n\in\N}\subset\Lambda_0$ such that $x_n\to y \in\Lambda_1$.
By corollary \ref{forma} $\Lambda_1$  is the disjoint union of  attractors and a repellers so suppose that $y\in\Lambda_1^s$.  For a sufficiently big $n$, $x_n$ must be in the basin of attraction  of $\Lambda_1^s$ and therefore $f^m(x_n)$ does not return to a neighborhood of $x_n$, which is impossible since  $ \Omega (f \mid_{\Lambda}) = \Lambda $.

\endproof

\vspace {0.2cm}
Now we are in the right conditions to prove:
\begin{theorem}
 Let $f: M\to M$ be a  diffeomorphism, $M$ a compact surface  and $\Lambda\subset M$  a  compact hyperbolic invariant set. If we also have that  $\Omega (f\mid_{\Lambda})=\Lambda$ then for any neighborhood  $V$ of $\Lambda$, there exist  $\widetilde{\Lambda}$ a  compact hyperbolic invariant set   with local product structure such that, $$\Lambda\subset\widetilde{\Lambda}\subset V\,.$$

\end{theorem}
\proof
Let $V$ be any open set containing $\Lambda$ and $V'=V\setminus\Lambda_1$. Note that $V'$ is open and contains $\Lambda_0$.

From lemma \ref{ultimo}, we have that there exists $\Lambda_0'$ with local product structure such that $\Lambda_0'\cap\Lambda_1=\emptyset$ and $ \Lambda_0 \subset \Lambda'_0 \subset V' $. On the other hand $\Lambda_1$ has local product structure from corollary \ref{cor1}.

We conclude that  $\widetilde{\Lambda}= \Lambda'_0\cup\Lambda_1 $ has local product structure, is contained in $V$ and contains $\Lambda$.
\endproof

\section*{Acknowledgment}
 This work is part of master thesis
at the Universidad de la Repbulica Uruguay under the guidance of Martín
Sambarino, I would not have been able to write this paper without his help support and infinite patience.
 I can not find a way to describe just how much he contributed to this work or how much he contributed to my education in general, so probably I should just say thanks for all. He also helped me greatly in the writing of the paper, dedicating a lot of work to correcting and re-correcting it. I would like to thank Rafael Potrie for repeatedly listening to me, and helping me understand. Also for suggesting ideas and encouraging me to work harder. To Javier Correa who helped with the writing of the paper. And finally to Todd Fisher and Keith Burns who kindly listened to the results and
showed interest in it, also they suggested me ideas that I am still thinking about in order to possibly understand these examples further.
\newpage

Cmat- Faculad  de Ciencias-Universidad de la República.

Igua 4225
Montevideo 11400
Uruguay

email: adaluz@cmat.edu.uy
\end{document}